\documentclass{amsart}

\usepackage{amssymb, bm, amsfonts, stmaryrd, amsmath, amsthm, enumerate, mathtools, mathrsfs, hyperref, colonequals, graphicx, color, nccmath}

\usepackage[capitalise]{cleveref}
\usepackage[parfill]{parskip}
\usepackage[margin=1.2in]{geometry}
\usepackage{multicol}
\usepackage{wrapfig,tikz,tikz-cd}
\usetikzlibrary{arrows}
\usepackage[foot]{amsaddr}

\usepackage[normalem]{ulem}

\theoremstyle{plain}\newtheorem{theorem}{Theorem}[section]
\newtheorem{theorema}{Theorem}
 
\theoremstyle{plain}\newtheorem*{conjecture}{Conjecture}
\theoremstyle{plain}\newtheorem{corollary}[theorem]{Corollary}
\theoremstyle{plain}\newtheorem{lemma}[theorem]{Lemma}
\theoremstyle{plain}\newtheorem{proposition}[theorem]{Proposition}
\theoremstyle{plain}
\theoremstyle{definition}\newtheorem{definition}[theorem]{Definition}
\theoremstyle{definition}
\theoremstyle{definition}
\theoremstyle{definition}\newtheorem{notation}[theorem]{Notation}
\theoremstyle{definition}\newtheorem{remark}[theorem]{Remark}
\theoremstyle{definition}
\theoremstyle{definition}\newtheorem{question}[theorem]{Question}

\DeclareMathOperator{\HH}{HH}

\DeclareMathOperator{\Ext}{Ext}

\DeclareMathOperator{\Der}{Der}
\DeclareMathOperator{\Inn}{Inn}
\DeclareMathOperator{\im}{im}

\renewcommand{\SS}{\mathcal{S}}
\newcommand{\WW}{\mathcal{W}}
\newcommand{\EE}{\mathcal{E}}
\newcommand{\ES}{\mathcal{ES}}

\title{Outer derivations on blocks of group algebras}

\author{Matthew Antrobus$^1$}
\address{
$^1$Department of Mathematics, University of Manchester, Manchester M13 9PL, United Kingdom {\tt matthew.antrobus@manchester.ac.uk}}

\author{Benjamin Briggs$^2$}
\address{
$^2$Department of Mathematics, Imperial College London, London SW7 2AZ, United Kingdom {\tt ben.briggs@imperial.ac.uk}}

\author{Lleonard Rubio y Degrassi$^3$}
\address{
$^3$Dipartimento di Matematica, Università degli Studi di Padova, Via Trieste 63, 35121 Padova PD, Italy {\tt lleonard@math.unipd.it} {\rm (corresponding author)}}

\thanks{\noindent
{\it Funding.}
The first author was funded by the Heilbronn Institute for Mathematical Research. The third author was funded by the INDAM- GNSAGA Project (CUP E53C24001950001) and also by the Department of Mathematics ‘Tullio Levi-Civita’ of the University of Padova through its SID 2025 - Investimento Strategico di Dipartimento, within the project `Homotopical and relative methods in Hochschild cohomology'.\\
{\it Acknowledgements.}
We would like to thank  Markus Linckelmann and Jay Taylor for their helpful comments and questions. 
}

\begin{document}

\begin{abstract}
Let $B$ be a block of a finite group algebra $kG$. Linckelmann has conjectured that if the defect group $P$ of $B$ is  nontrivial then $B$ admits a $k$-linear outer derivation, or in other words $\HH^1(B)\neq 0$.

We provide various criteria for the non-vanishing of $\HH^1(B)$ in terms of the subgroup $P\subseteq G$. In particular, we show that Linckelmann's conjecture holds for principal blocks having abelian defect group, for all blocks of the symmetric and  alternating groups, for blocks of finite groups of Lie type in defining characteristic, and for blocks of general linear groups in any characteristic. Our main tool relates the non-vanishing of the Batalin-Vilkovisky operator  $\Delta\colon \HH^1(B)\to \HH^0(B)$ to the existence of a certain elements of $P$, called extra-strong non-Schur elements.

We show that extra-strong non-Schur $p$-elements exist for all groups having order divisible by a prime $p>5$, extending work of Fleischmann, Janiszczak, and Lempken. Using this, we prove that if $k$ has characteristic $p>5$, then Linckelmann's conjecture holds for all blocks of $kG$ with Sylow defect group. Likewise, when $k$ has characteristic $p>5$ we deduce that the first Hochschild cohomology is non-zero for any nontrivial twisted group algebra over $k$.

\end{abstract}
\raggedbottom

\maketitle

\maketitle

\section{Introduction}

Determining the existence, or non-existence, of outer derivations on algebraic structures has led to important work in a number of areas. Whitehead's First Lemma for derivations on Lie algebras \cite[Chapter III, \S 7, Lemma 3]{Jacobson62}, or the derivation problem for locally compact groups, eventually solved in \cite{Losert}, are notable examples. Derivations on finite group algebras in particular have been an active research topic since the 1970s \cite{Smith, Spiegel,AMS, AM}, and, for example, have seen applications to coding theory in computer science \cite{Bu,CH}.

Let $G$ be a finite group and let $k$ be a field whose characteristic divides the order of $G$. In this context the derivation problem asks whether the space $\HH^1(kG) = \Der(kG,kG)/\Inn(kG,kG)$, consisting of $k$-linear derivations modulo inner derivations, is nonzero. It was proven by Gerstenhaber and Green that if $kG$ admits a semisimple deformation, that is, if the  Donald-Flanigan conjecture holds for $G$, then $\HH^1(kG)\neq 0$ \cite{GerstenhaberGreen}. While the Donald-Flanigan conjecture remains open, we know that $\HH^1(kG)$ is always nonzero thanks to work of Fleischmann, Janiszczak, and Lempken; this is a consequence of a stronger, group-theoretic statement proven in \cite{FJL} using the classification of finite simple groups.

Perhaps in search of a more conceptual argument for the existence of outer derivations,  Linckelmann posed the block-wise version of this result as a question in {\cite[Question 7.4]{Linckelmann2}}, and later stated it as a conjecture in \cite{LinckelmannTalk}:

\begin{conjecture}[Linckelmann]
\label{LConj}
If $B$ is a block of a finite group algebra $kG$ with a nontrivial defect group then  $\mathrm{HH}^1(B)\neq 0$.
\end{conjecture}

A number of authors have shown that this conjecture holds in special cases. For instance if the Gabriel quiver of $B$ has a loop \cite[Theorem 1.2]{BKL}; if $B$ is block of a symmetric group algebra $kS_n$ with an abelian defect group \cite[Example 7.5]{Linckelmann2}; and more generally, if $B$ is any block of $kS_n$ \cite[Theorem 1.1]{BKL23}; for blocks of the sporadic Mathieu groups \cite[Theorem 1.1]{Murphy2}; for a direct sum of blocks of some central group extension \cite{Murphy}; for blocks with normal defect group \cite[Corollary 1.4]{Todea}; and for most blocks of the Janko groups \cite[Theorem 1.6]{ST}.  In addition, in \cite[Theorem 1.5.]{ST} the authors obtain a partial reduction theorem to quasi-simple groups.

Denote by $E^p(G)$ the smallest normal subgroup of $G$ such that $G/E^p(G)$ is an elementary abelian $p$-group. An element of $G$ is a \emph{extra non-Schur element} if $x\notin E^p(C_G(x))$.  Let $P$ be a $p$-subgroup of $G$. If $x\in P$ we say that $G$ has the \emph{extra-strong non-Schur property relative to} $P$. We write $\ES(G, P)$ for the set of $G$-conjugacy classes of such elements. 

We approach Linckelmann's problem by utilising the Batalin-Vilkovisky operator 
\[
\Delta_{A}^1\colon \HH^{1}(A) \rightarrow \HH^0(A),
\]
defined originally for symmetric algebras by Tradler in \cite{Tradler}, and studied in the context of group algebras by Benson, Kessar and Linckelmann in \cite{BKL21}. Given a block $B$ of $kG$ with defect $P$, we then relate the non-vanishing of $\Delta_B^1$ to the existence of extra-strong non-Schur elements in $G$ relative to $P$.
\begin{theorema}\label{thm_1_cor}
    Let $B$ be a block of $kG$ having defect group $P$. Then $\operatorname{dim}(\HH^1(B)) \geq |\ES(G,P)|$.
\end{theorema}

We use this result to establish many more cases in which $\HH^1(B)$ is nonzero:

\begin{theorema}\label{thm_2}
    Let $k$ be a field of characteristic $p$, let $G$ be a finite group, and let $B$ be a block of $kG$ with a nontrivial defect group $P$. In all of the following cases $\HH^1(B)\neq 0$:
    \begin{enumerate}
    \item 
    The map ${\rm H}_1(P,\mathbb{F}_p)\to {\rm H}_1(G,\mathbb{F}_p)$ is nonzero.
    \item  
    There is $x\in P\smallsetminus E^p(P)$ such that $C_P(x)$ is a Sylow subgroup of $ C_G(x)$.
        \item \label{abelian} $P$ is Sylow (for example if $B$ is the principal block) and abelian.
        \item \label{item_nonSchurelt1}
        $P$ is Sylow, and $G$ contains a extra  non-Schur $p$-element.
    \item \label{item_an}  $G$ is a symmetric group $S_n$ or alternating group $A_n$.
    \item  
    $G$ is a finite general linear group $\mathrm{GL}_n(\mathbb{F}_q)$, defined in any characteristic $r$ with $q=r^n$. 
    \item 
    $G$ is a finite group of Lie type ${\bf G}^F$, or the commutator subgroup $[{\bf G}^F, {\bf G}^F]$, where ${\bf G}$ is a reductive algebraic group over $\overline{\mathbb{F}}_p$, split over $\mathbb{F}_p$,  
    and where $F$ is a Frobenius map.
    \item
    There is a normal $p$-subgroup $Q\subseteq G$ inducing a nonzero map ${\rm H}_1(Q,\mathbb{F}_p)\to {\rm H}_1(G,\mathbb{F}_p)$.
    \end{enumerate}
    Note that the last four parts only depend on $G$, not on the block or its defect group.
\end{theorema}

See \cref{prop_nonSchur,thm:HH1blockex,thm_Sncase,thm_Ancase} for these results (and several more). We also recover a number of results from the literature; see \cref{rmk:knowncases,thm_Sncase}. 
Brou\'e's abelian defect group conjecture states that  any block $B$ of $kG$ with abelian defect group $P$ is derived equivalent to its Brauer correspondent, which is a block of $kN_G(D)$  \cite[Question~6.2]{Broue}. The validity of this conjecture would imply that $\HH^1(B)\neq 0$ whenever $B$ has abelian defect group (see Remark \ref{rmk:knowncases}) so we see (\ref{abelian}) as evidence for Brou\'e's abelian defect group conjecture.

Fleischmann, Janiszczak, and Lempken in \cite{FJL} proved that many (but not all) finite groups possess non-Schur $p$-elements,  that is, a $p$-element $x$ such that $x \notin [C_G(x),C_G(x)]$.  Extending their computations we verify, using the classification of finite simple groups, that all finite groups possess extra-strong non-Schur $p$-elements as long as $p>5$. As a consequence, in Section \ref{sec_char_5} we establish Linckelmann's conjecture for all  blocks of maximal defect in characteristic larger than $5$.

\begin{theorema}\label{thm_3}
If $p$ is a prime number greater than $5$ then any finite group $G$ with order divisible by $p$ contains a extra-strong non-Schur $p$-element. Consequently, if $k$ is a field of characteristic $p$, then any block $B$ of $kG$ of maximal defect (e.g.\ the
principal block) satisfies $\HH^1(B)\neq 0$.
\end{theorema}

In \cref{sec_Hochschild} we give a different technique  for establishing $\HH^1(B)\neq 0$. The subgroup inclusion $P\subseteq G$ induces maps of Hochschild cohomology groups
\[
\HH^1(kP,kP)\xrightarrow{\alpha} \HH^1(kP,kG)\xleftarrow{\beta} \HH^1(kG,kG),
\]
given by post- and precomposing a derivation with the algebra homomorphism $kP\to kG$. 

We prove, as consequence of some general results about the Hochschild cohomology groups and separable equivalences,  that for any block $B$ of $kG$ with defect group $P$, there is an embedding of $\alpha^{-1}(\im(\beta))$ into $\HH^1(B)$. In other words, we prove 
\begin{theorema}\label{thm_alpha_beta}
    The space of derivations on $kP$ which are restrictions of derivations on $kG$, modulo inner derivations, is isomorphic to a subspace of $\HH^1(B)$. 
\end{theorema}
 
 Denote by $\mathrm{Core}_G(P)$ the largest normal subgroup of $G$ contained in $P$ and by $O_p(G)$ of $G$ the largest normal $p$-subgroup of $G$. From \cref{thm_alpha_beta} we deduce

\begin{theorema}
Let $G$ be a finite group and let $B$ be a block of $kG$ with defect group $P$.  In all of the following cases $\HH^1(B)$ is non-zero:
\begin{enumerate}
    \item If for some $x$ in $\mathrm{Core}_G(P)$, the homomorphism ${\rm H}_1(C_P(x),\mathbb{F}_p)\to {\rm H}_1(C_G(x),\mathbb{F}_p)$ is nonzero.
    \item  If there is a $p$-element $x$ in $O_p(G)$ such that ${\rm H}_1(C_G(x),\mathbb{F}_p)\neq 0$ and  $B$ is a block  of maximal defect group.
    \item If $O_p(G)\cap Z(P)$ is not contained in $[P,P]$ and  $B$ is a block of maximal defect group. 
\end{enumerate}
\end{theorema}

In the final part of the paper, Section \ref{s:tga}, we discuss twisted group algebras. The relevance of these to Linckelmann's Conjecture lies in the fact in many cases, for instance when  $B$ has normal defect group or if $G$ is p-solvable,  blocks of $kG$ are Morita equivalent to twisted group algebras.  
Based on this, Todea posed the following question:  
   ``\textit{Is the first Hochschild cohomology of twisted group algebras of finite groups
(with order divisible by p) nontrivial? If not, can someone find examples of such twisted
group algebras which have zero first Hochschild cohomology?}" \cite[Question 1.3]{Todea}. A partial affirmative answer to this question has been obtained  in various cases, such as $p$-solvable
groups \cite[Corollary 1.4]{Todea}, certain dihedral groups \cite{OteroSanchez25},  and all finite simple groups \cite[Corollary 1.2]{Murphy} (note that no reduction to finite simple groups is known). As a consequence of \cref{thm_3} we obtain:
\begin{theorema}\label{thm_4}
    Let $G$ be a finite group and let the characteristic of the field $k$ be greater than $5$. Then  $\mathrm{HH}^1(k_\alpha G,k_\alpha G) \neq 0$ for every $2$-cocycle $\alpha \in Z^2(G; k^\times)$.
\end{theorema}

\section{Background}
In this section we collect together some group theoretic notations that will be useful later, as well as some facts about Hochschild cohomology. 

Throughout, $G$ is a finite group and  $k$ is a field of positive characteristic $p$. All algebras are $k$-algebras, and all unadorned tensor products and cohomology groups are taken relative to $k$.

The cohomology of $G$ with coefficients in $k$ is $\operatorname{H}^*(G,k)$, and we will freely make the identifications $\operatorname{H}^1(G,k)= \operatorname{Hom}(G,k^+)$ and $\operatorname{H}^0(G,k) = k$.

The subgroup $E^p(G):=G^p[G,G]$ is the smallest normal subgroup of $G$ such that $G/E^p(G)$ is an elementary abelian $p$-group. For any element $x\in G$ we note that $x\notin E^p(G)$ if and only if  $[x]\neq 0$ in ${\rm H}_1(G,\mathbb{F}_p)$.

Let $x$ be an element of $G$. The centralizer of $x$ in $G$ is written $C_G(x)$. The conjugacy class of $x$ is $C_x$. Given $U \subseteq G$, we will write $\widehat{U} := \sum_{g \in U} g$. So for example, $\widehat{C}_x$ is the sum of all elements conjugate to $x$. The element $x$ is a $p$-element if its order if a $p$-power, and $x$ is {$p$-regular} if its order is not divisible by $p$. Any $x\in G$ can be written uniquely as $x=x_px_{p'}$, where $x_p$ is a $p$-element and $x_{p'}$ is $p$-regular. The $p$-section of $x$ in $G$ consists of those elements $y \in G$ for which $y_p$ is conjugate in $G$ to $x_p$.

For any $U \subseteq G$ we write $\pi_U\colon kG \rightarrow kG$ for the projection $\pi_U(a) = \sum_{g \in U} \alpha_g g$, where $a = \sum_{g \in G} \alpha_g g$. The Brauer map $\operatorname{Br}_P\colon Z(kG) \rightarrow Z(kC_G(P))$ associated to a $p$-subgroup $P\subseteq G$ is the restriction of $\pi_{C_G(P)}\colon kG\to kC_G(P)$ to the centre $Z(kG)$.

The group algebra of $G$ decomposes uniquely into a product $kG=B_1\times\cdots \times B_n$ of connected rings $B_i$, called the blocks of $kG$. Each block $B_i$ corresponds uniquely to a central idempotent $e_{B_i}$ in $kG$. A defect group for $B_i$ is a maximal $p$-subgroup $D\subseteq G$ for which $\operatorname{Br}_D(e_{B_i})\neq 0$.

We write $EG \subseteq Z(kG)$ for the subring of $Z(kG)$ generated by central idempotents in $G$. Equivalently, $EG$ is the linear span of $e_{B_1},\ldots,e_{B_n}$.

\begin{lemma}[{\cite[(61)]{Kulshammer91}}]\label{Kulshammer91}  If $U$ is a $p$-section in $G$ then the restriction of $\pi_U$ to $Z(kG)$ is an $EG$-module homomorphism. 
\end{lemma}

\begin{proposition}\label{proposition:ConjSumLindep} Let $B$ be a block of $kG$ having defect group $D$, and block idempotent $e_B$. Then the set $\{ \widehat{C}_x e_B~:~ x \in D\}$ is linearly independent in $Z(B)$. 
\end{proposition}
\begin{proof}
    Write $z = \sum_{ x\in D}  \lambda_x\widehat{C}_x$ and assume that $ze_B= 0$. Our goal is to show that each $\lambda_x = 0$.
    
    Write $\operatorname{Br}_{\langle x \rangle}\colon Z(kG) \rightarrow Z(kC_G(x))$ for the Brauer map. Note that $\operatorname{Supp}(\operatorname{Br}_{\langle x \rangle}(z))\subseteq\operatorname{Supp}( z) \subseteq \bigcup_{y\in D} C_y$ consists entirely of $p$-elements.

    Now let $U$ be the $p$-section of $x$ in $C_G(x)$. If $y\in U$ is a $p$-element then $y$ is conjugate to $x$ in $C_G(x)$, but $x$ is central in $C_G(x)$, so $y=x$. Hence $\operatorname{Supp}(\operatorname{Br}_{\langle x \rangle}(z)) \cap U\subseteq \{x\}$, and so $\pi_U(\operatorname{Br}_{\langle x \rangle}(z))=\lambda_x x$.

 Since the Brauer map is a ring homomorphism, $\operatorname{Br}_{\langle x \rangle}(e_B)$ is a central idempotent of $kC_G(x)$. From here, using the fact that $\pi_U$ is an $EC_G(x)$-module homomorphism by \cref{Kulshammer91}, we have
 \[
 0=\pi_U(\operatorname{Br}_{\langle x \rangle}(ze_B)) = \pi_U(\operatorname{Br}_{\langle x \rangle}(z)\operatorname{Br}_{\langle x \rangle}(e_B)) = \pi_U(\operatorname{Br}_{\langle x \rangle}(z))\operatorname{Br}_{\langle x \rangle}(e_B) = \lambda_x x\operatorname{Br}_{\langle x \rangle}(e_B).
 \]
From \cite[Theorem 4.11]{Navarro98} $\operatorname{Br}_{\langle x \rangle}(e_B)$ is non-zero, and so this implies $\lambda_x=0$.
\end{proof}

\subsection{Hochschild cohomology}
If $A$ is an associative algebra and $M$ is an $A$-bimodule, then the Hochschild cohomology of $A$ with coefficients in $M$ is 
\[
\HH^*(A,M)=\Ext^*_{A^{\rm op}\otimes A}(A,M).
\]
We will frequently use the short-hand notation $\HH^*(A)=\HH^*(A,A)$, and we will freely make the identifications $\HH^0(A) = Z(A)$ and $\HH^1(A) = \Der_k(A,A)/\Inn(A,A)$, 

We can concretely describe the Hochschild cohomology of a group algebra in the following way.

\begin{proposition}[Centraliser Decomposition]\label{prop_centraliser_decomposition} We have an isomorphism of vector spaces
\[
\HH^*(kG) \cong \bigoplus_{[x]\in G/ \sim_G} \operatorname{H}^*(C_G(x), k).
\]
If we write $\gamma_x^*$ for the inclusion $\operatorname{H}^*(C_G(x),k) \hookrightarrow \HH^*(kG)$ induced by this decomposition, then $\gamma_x^0$ sends $\lambda \in k$ to $\lambda \widehat{C}_x\in Z(kG)$.
\end{proposition}

\begin{proof} The isomorphism of vector spaces is proved in \cite[Theorem 2.11.2]{Benson_1991}. We verify the second claim. Given a $kG$-bimodule $M$, write $M^{\rm ad}$ for the $kG$-module $M$ with the  conjugation action of $G$. We have $\HH^0(kG) \cong \operatorname{H}^0(G,kG^{\rm ad})$, and $kG^{\rm ad} \cong \bigoplus_{x \in G/\sim} kC_x^{\rm ad} \cong \operatorname{Ind}_{C_G(x)}^Gk$. So we see that
$$k \xrightarrow{=} \operatorname{H}^0(G,kC_x^{\rm ad}) \rightarrow Z(kG)$$
is (up to a choice of scalar) the map we claimed.
\end{proof}
Theorem \ref{thm_1_cor}  will be derived from non-vanishing of a certain map on Hochschild cohomolology, called the Batalin-Vilkovisky (BV) operator.
\begin{theorem}\cite[Theorem 1]{Tradler} Let $A$ be an algebra equipped with a symmetric, invariant, non-degenerate bilinear form $\langle -, -\rangle_A$. Then there exists an operator $\Delta_A\colon \HH^n(A) \rightarrow \HH^{n-1}(A)$, endowing $\HH^*(A)$ with the structure of a Batalin-Vilkovisky algebra.  \end{theorem}
When $A$ is a group algebra $kG$, some properties of the BV operator were studied in \cite{BKL21}, where in particular the following was shown. 
\begin{theorem}\cite[Theorem 1.1]{BKL21} \label{theorem:BVBKL} The BV operator on the Hochschild cohomology of $kG$ preserves the centralizer decomposition, and restricts to maps 
$$\Delta_x\colon \operatorname{H}^i(C_G(x),k) \rightarrow \operatorname{H}^{i-1}(C_G(x),k),$$
so that $\gamma_x \Delta_x = \Delta_{kG} \gamma_x$. Furthermore, given $\psi \in \operatorname{H}^1(C_G(x),k)=\operatorname{Hom}(C_G(x),k)$, $\Delta_x(\psi) = \psi(x)$. 
\end{theorem}
With careful choice of inner products, the BV operator respects the block decomposition of Hochschild cohomology.
\begin{proposition}\label{proposition:BlockDecompBV} Given a block $B$ of $kG$ with block idempotent $e_B$, and associating the natural inner product to $B$, the following diagram commutes
\begin{center}
    \begin{tikzcd}
        \HH^1(kG) \arrow[r, "\Delta_{kG}"] \arrow[d, "-\cdot e_B"] & \HH^{0}(kG) \arrow[d, "-\cdot e_B"] \\
        \HH^1(B) \arrow[r, "\Delta_B"] & \HH^{0}(B).
    \end{tikzcd}
\end{center}

\end{proposition}
\begin{proof}
Given $f \in \operatorname{HH}^1(A)$, the BV operator can be defined as the element satisfying $\langle \Delta_A(f), x \rangle_A = \langle f(x), 1\rangle_A$. Thus given $f \in \HH^1(kG)$ and $x \in B$, we have
\[
\langle \Delta_{kG}(f)e_B, x \rangle_{kG} =\langle \Delta_{kG}(f),e_B x \rangle_{kG}= \langle f(e_Bx), 1 \rangle_{kG} = \langle e_Bf(x), e_B\rangle_B = \langle \Delta_B(e_Bf), x \rangle_B,
\]
where we have used that derivations are $EG$-homomorphisms. By non-degeneracy of the inner product, we are done.
\end{proof}

\section{Criteria for non-vanishing of Hochschild cohomology}

In this section we establish some group theoretic criteria for the non-vanishing of the first Hochschild cohomology of a block. These criteria are based on variations of the non-Schur property considered by Fleischmann, Janiszczak, Lempken in \cite{FJL}; the next definition summarises all of these variations.

\begin{definition}
    Let $G$ be a finite group, let $p$ be a prime number, and let $x\in G$.
    
  Following \cite{FJL}, the element $x$ is a \emph{non-Schur element} of $G$ if $x \notin [C_G(x),C_G(x)]$. Further, $x$ is a \emph{weak non-Schur element} if it is  non-Schur with order divisible by $p$,   while $x$ is a \emph{strong non-Schur element} if it is a non-Schur $p$-element. We write $\WW(G,p)$ and $\SS(G,p)$ for the set of conjugacy classes of weak and strong non-Schur elements correspondingly.

  Extending this terminology, we say that $x$ is a \emph{extra non-Schur element} if $x$ is an element 
  for which $x \notin E^p(C_G(x))=C_G(x)^p[C_G(x),C_G(x)]$, or equivalently, if $[x]\neq 0$ in ${\rm H}_1(C_G(x),\mathbb{F}_p)$. As before, $x$ is an \emph{extra-strong non-Schur element} if it is also a $p$-element.   We write $\EE(G,p)$ and $\ES(G,p)$ for the set of conjugacy classes of extra non-Schur and extrastrong non-Schur elements correspondingly.

    We say that $G$ has the \emph{weak non-Schur property} $W(p)$ if it contains a weak non-Schur element,  or if $p\nmid |G|$ (we consider $W(p)$ to hold vacuously in the latter case). Likewise $G$ has the \emph{strong non-Schur property} $S(p)$ or the \emph{extra-strong non-Schur property} $ES(p)$ 
    when it possess an element of the corresponding kind.
    
    More generally, given a $p$-subgroup $P$ of $G$, we say that $G$ has the \emph{extra non-Schur property relative to} $P$ if there is extra non-Schur $x \in G$ with $x_p \in P$. Similarly, we say that $G$ has the \emph{extra-strong non-Schur property relative to} $P$ when there is an extra-strong $x \in P$. We write $\mathcal{E}(G,P)$ for the set of $G$-conjugacy classes of extra non-Schur elements relative to $P$, and $\ES(G,P)$ for the set of $G$ conjugacy classes of extra-strong non-Schur elements relative to $P$. Note that $\ES(G,P)=\ES(G,p)$ when $P$ is a Sylow subgroup of $G$.

\end{definition}

\begin{theorem}\label{thm:main1'}
    Let $B$ be a block of $kG$ having defect group $P$. Then $\operatorname{dim}(\HH^1(B)) \geq |\ES(G,P)|$.
\end{theorem}
\begin{proof} Making use of the BV operator $\Delta_B^1\colon \HH^1(B)\to  \HH^{0}(B)$, we establish the claimed inequality by demonstrating that $\operatorname{dim}(
\operatorname{Im}(\Delta_B^1)) \geq |\ES(G,P)|$.

Using \cref{theorem:BVBKL,prop_centraliser_decomposition} we know that $\Delta_{kG}(\gamma_x(\psi)) = \psi(x)\widehat{C}_x$ for any conjugacy class $C_x$ and any $\psi\in \operatorname{H}^1(C_G(x), k)$. So then, by \cref{proposition:BlockDecompBV}, we have $\Delta_B(\gamma_x(\psi)e_B) = \Delta_{kG}(\gamma_x (\psi)) e_B = \psi(x)\widehat{C}_x e_B $. If $C_x\in \ES(G,P)$ then there is some $\psi \in \operatorname{H}^1(C_G(x),k)$ such that $\psi(x)\neq 0$, and therefore $\widehat{C}_xe_B \in \operatorname{Im}(\Delta_B^1)$.  By \cref{proposition:BlockDecompBV} as we vary $C_x \in \ES(G,P)$ the elements $\widehat{C}_xe_B$ are linearly independent. 
\end{proof}

\begin{remark}
   Our theorem restricts attention to extra-strong non-Schur elements, though non-vanishing of outer derivations could be determined by extra-non-Schur elements. Indeed, $J_4$ at $p=3$ admits no extra-strong non-Schur elements, but the presence of extra non-Schur elements demonstrates non-vanishing of the blocks of maximal defect. By \cite[Proposition 7]{Kulshammer94}, if $x \in G$ has $p$-part not conjugate to an element of $P$ then $e_B \widehat{C}_x=0$. 
   We then  have 
   \[
   \operatorname{span}_k\big\{e_B\widehat{C} ~:~ C\in \ES(G,p )\big\} \subseteq  \operatorname{span}_k\big\{e_B\widehat{C} ~:~ C\in \mathcal{E}(G,p )\big\}=\operatorname{Im}(\Delta^1_B)
   \]
   Further \cite[Theorem 6]{Kulshammer94} implies that $\widehat{C}_x e_B$ is determined by $\operatorname{Br}_{\langle x \rangle}(\widehat{C}_x)$ and an idempotent of $kC_G(x)$. 
\end{remark}
\begin{remark} The presence and quantity of extra-strong non-Schur elements in $G$ can be deduced from the fusion system $\mathcal{F}=\mathcal{F}_P(G)$ of $G$ on its Sylow subgroup $P$, in the following manner. Conjugacy classes of $p$-elements of $G$ give rise to $\mathcal{F}$-isomorphism classes of fully centralized cyclic subgroups in $\mathcal{F}_P(G)$. One can check that given such a subgroup $\langle x\rangle$ we have $C_{\mathcal{F}}(\langle x \rangle) = \mathcal{F}_{C_P(x)}(C_G(x))$. Then, noting that $\operatorname{H}^*(C_G(x),\mathbb{F}_p) \cong \operatorname{H}^*(C_{\mathcal{F}}(\langle x \rangle),\mathbb{F}_p) \cong \varprojlim_{C_\mathcal{F}(\langle x \rangle)} \operatorname{H}^*(-, \mathbb{F}_p)$, the conjugacy class of $x \in G$ is extra-strong non-Schur exactly when $\operatorname{H}^1(C_\mathcal{F}(\langle x \rangle), k) \rightarrow \operatorname{H}^1(\langle x \rangle, \mathbb{F}_p)$ is non-zero, where this map is induced by the inverse limit. Each such isomorphism class of subgroup contributes $(p-1)/{|\operatorname{Aut}_\mathcal{F}(\langle x \rangle)|}$ extra-strong non-Schur elements.
\end{remark}

We record two immediate consequences of \cref{thm:main1'}.

\begin{corollary}\label{thm:main1}
        Let $B$ be a block of $kG$ with defect group $P$. If there is $x\in P$ such that $x\notin E^p(C_G(x))$ then $\HH^1(B)\neq 0$.
\end{corollary}

\begin{corollary}\label{cor_nonzero_Sylow}
    Let $G$ be a finite group. If there is a $p$-element $x\in G$ such that $x\notin E^p(C_G(x))$ then any block $B$ of $kG$ of maximal defect group satisfies $\HH^1(B)\neq 0$.
\end{corollary}

 \begin{remark}\label{rmk:Scp}
If $G$ possesses  an extra-strong non-Schur element then by \cref{cor_nonzero_Sylow} any block $B$ of $kG$ of maximal defect group satisfies $\HH^1(B)\neq 0$. The converse does not hold. For example, if $G$ is the  Thompson group ${\rm Th}$ and $k$ is a field of characteristic $p=5$, then the principal block $B$ satisfies $\HH^1(kG)=\HH^1(B)\neq 0$ although $G$ fails the strong commutator index property \cite[P3]{Todea}. In particular $G$ does not satisfy $\ES(p)$.
\end{remark}

The following (well-known) lemma will be useful later.

\begin{lemma}\label{lem:CPx}
    Let $H\subseteq G$ be a subgroup with index coprime to $p$ and assume that $x\in H$ is central in $G$. Then $x\in E^p(H)$ if and only if $x\in E^p(G)$.
   \end{lemma}

\begin{proof}
    Since $E^p(H)\subseteq E^p(G)$ it is clear that $x\in E^p(H)$ implies $x\in E^p(G)$.

     We recall that the transfer  map $\operatorname{tr}^{G}_H\colon G\to H^{\rm ab}$ is a group homomorphism defined by choosing right coset representatives $\{\gamma_1, \dots, \gamma_n\}$ for $H$ in $G$; there are $x_i\in H$ for which $\gamma_ix=x_i\gamma_j$ for some $j$, and by definition
    \[
    \operatorname{tr}^{G}_H(x)= x_1\cdots x_n.
    \]
    In our particular case $x$ is central in $G$ and we have $\gamma_ix=x\gamma_i$, so $ \operatorname{tr}^{G}_H(x)= x^n.$

    If  $x\notin E^p(H)$ then there is a group homomorphism $\varphi:H\to k^+$ such that $\varphi(x)\neq 0$. The composition $\varphi\circ\operatorname{tr}^{G}_H:G \to k^+$ maps $x$ to $n\varphi(x)\neq 0$, since $n=[G:H]$ is coprime to $p$. Hence $x\notin E^p(G)$.
\end{proof}

We can give another criterion which only depends on the fusion system $\mathcal{F}_P(G)$ of $G$ on $P$.

\begin{corollary}
   Let $G$ be a finite group with Sylow subgroup $P$. Assume there is a minimal generator $x\in P$ for which the subgroup $\langle x \rangle$ generated by $x$ is  fully centralised in $\mathcal{F}_P(G)$. Then any block $B$ with maximal defect satisfies $\HH^1(B)\neq 0$.

\end{corollary}

\begin{proof}
    By \cite[Proposition 1.3.]{BLO03} the subgroup $\langle x \rangle$ is fully centralised  in $\mathcal{F}_P(G)$ if and only if $C_P(\langle x \rangle)=C_P(x)$  is a Sylow $p$-subgroup of $C_G(x)$. Since $x$ is a minimal generator of $P$, it is also a minimal generator of $C_P(x)$, and so $x\notin E^p(C_P(x))$.     The statement now follows by combining \cref{cor_nonzero_Sylow} with \cref{lem:CPx}.
\end{proof}

\begin{definition}
    Following \cite{FJL} a finite group $G$ is said to have the \emph{weak non-Schur property} $W(p)$ if it possesses an weak non-Schur element, or if $p\nmid |G|$ (we consider $W(p)$ to hold vacuously in the later case). Likewise $G$ has the  \emph{strong} or  \emph{extra-strong non-Schur property} $S(p)$ or $ES(p)$ if it possesses an element of the corresponding kind, or if  $p\nmid |G|$.
\end{definition}

 In this work we will be most interested in the extra-strong non-Schur property $ES(p)$. As in the strong case in \cite[Proposition 1.3]{FJL}, the next proposition will imply that every minimal counterexample to $ES(p)$ is a simple group.

\begin{proposition}\label{Prop:red}
  If both $N$ and $G/N$ satisfy $ES(p)$, then so does $G$. In particular any minimal finite group not satisfying $ES(p)$ must be a simple group. 
\end{proposition}
\begin{proof}
     If $p$ does not divide the order of $N$ or $G/N$, then $ES(p)$ holds  vacuously for $G$.
     
     Assume $p$ divides $|G/N|$. Then by assumption $G/N$ admits  a extra-strong non-Schur element $yN$. Each Sylow $p$-subgroup of $G/N$ is the image of a Sylow $p$-subgroup of $G$, so as $yN$ lies in a Sylow $p$-subgroup, it is the image of a $p$-element $y \in G$. Since $yN\notin E^p(C_{G/N}(yN))$, there is a group homomorphism $C_{G/N}(yN) \to \mathbb{F}_p^+$ that does not vanish on $yN$. Hence $C_G(y)\to C_{G/N}(yN)\to \mathbb{F}_p^+$ does not vanish on $y$. Therefore $y\notin E^p(C_G(y))$.

Assume $p$ divides $|N|$ but not $|G/N|$. Then by assumption $N$ admits  a extra-strong non-Schur element $x$, which means $x \not \in E^p(C_N(x))$. Note that $C_G(x)/C_N(x)$ is a subgroup of $G/N$, and so the index $[C_G(x):C_N(x)]$ is coprime to $p$. By \cref{lem:CPx}, $x \not \in E^p(C_N(x))$ if and only if $x \not\in E^p(C_G(x))$.
\end{proof}

   The previous proposition means that a group satisfies $ES(p)$ as long as all of its simple composition factors do. The next proposition gives group theoretic criteria for $ES(p)$ to hold; most are slight generalizations of results from \cite{FJL} but we have added a few that will be useful to us later.

\begin{proposition}\label{prop_nonSchurRelative}
Let $G$ be a finite group with a Sylow $p$-subgroup $P$ with a nontrivial subgroup $D \subseteq P$. In the following cases, $D$ contains an extra-strong non-Schur element.

\begin{enumerate}[\hspace{5mm} (i)]
\item\label{item_PprimeRel} There is an element $x \in D$ and a Sylow subgroup $Q\subseteq C_G(x)$ such that $x\in Q\smallsetminus E^p(Q)$.  
    \item\label{item_centreRel} $C_D(P)$ is not contained in $E^p(P)$.
    \item \label{item_sylabRel} $P$ is abelian.
     \item \label{item_syl_expRel} The exponent of $[P,P]$ is less than the exponent of $D$.
\end{enumerate}

\end{proposition}
\begin{proof}
    (\ref{item_PprimeRel})  The statement follows from \cref{lem:CPx}. 

(\ref{item_centreRel}) Pick $x\in C_D(P)\smallsetminus E^p(P)$. Since $P$ is a  Sylow $p$-subgroup of $G$ and $x$ is central, $P$ is also a Sylow $p$-subgroup of $C_G(x)$.  It follows from \cref{lem:CPx} that $x\notin E^p(C_G(x))$.

(\ref{item_sylabRel})   
This follows from (\ref{item_centreRel}).

(\ref{item_syl_expRel}) Take an element $x\in D$ of maximal order, write $C:= C_P(x)$ and suppose $x\in E^p(C)$. Working in $C/[C,C]$, we must have $y \in C$ so that $x[C,C] = y^p[C,C]$. Equivalently, there is $y \in C$ and $z \in [C,C]$ so that $xy^{-p}=z$. Since $z\in [P,P]$ the order $p^\ell$ of $z$ is strictly less than the order of $x$. As $x$ and $y$ commute we must have $x^{p^\ell} = y^{p^{\ell+1}}\neq 1$. But then the order of $y$ is greater than the order of $x$, a contraction.
\end{proof}
\begin{proposition}\label{prop_nonSchur}
    Let $G$ be a finite group with a nontrivial Sylow $p$-subgroup $P$. In the following cases $G$ satisfies $\ES(p)$.
    \begin{enumerate}[\hspace{5mm} (i)]
    \item\label{item_centre} $P$ has a central minimal generator.
    \item \label{item_syl_meta} $P$ is metacyclic. 
    \item \label{item_psolvable}$G$ is $p$-solvable.        
    \item \label{item_normal_com}
    $[P,P]$ is contained in a normal $p$-subgroup of $G$ (this holds in particular if $P$ is normal). 
    \end{enumerate}
\end{proposition}
\begin{proof}

(\ref{item_centre}) This follows by taking $P=D$ in \cref{prop_nonSchurRelative} (\ref{item_centreRel}), since minimal generators of $P$ are the same thing as elements of $P\smallsetminus E^p(P)$.

 (\ref{item_syl_meta}) Since $P$ is metacyclic, it has a cyclic normal subgroup $N$ such that $P/N$ is cyclic.  We want to show that the exponent $\mathrm{exp}([P,P])$ is less than $\mathrm{exp}(P)$ so the statement would follow from  (\ref{item_syl_expRel}).
 
 Assume by contradiction  $\mathrm{exp}([P,P])=\mathrm{exp}(P)$. Since $P/N$ is abelian, $[P,P]$ is a subgroup of $N$.  Since $\mathrm{exp}([P,P])\leq \mathrm{exp}(N)\leq \mathrm{exp}(P)$. Our assumption implies $\mathrm{exp}([P,P])=\mathrm{exp}(N)$, hence $[P,P]=N$. Since $P/[P,P]\subseteq P/N$ is cyclic, by the Burnside basis theorem $P$ is cyclic as well. Hence $\mathrm{exp}([P,P])<\mathrm{exp}(P)$, a contradiction.

 (\ref{item_psolvable}) The statements follow from \cref{Prop:red}.
 
(\ref{item_normal_com}) If $N\subseteq G$ is a normal $p$-subgroup containing $[P,P]$ then it satisfies $ES(p)$ by virtue of being a $p$-group. Moreover $G/N$ has abelian Sylow subgroups, so it satisfies $ES(p)$ according to (\ref{item_sylabRel}). Therefore $G$ satisfies $ES(p)$ by \cref{Prop:red}. 
\end{proof}

By definition $ES(p)$ implies $S(p)$. The following lemma provides some conditions for the converse to hold.

\begin{lemma}\label{lemma:StrongSuperStrong} Let $p$ be a prime, $G$ a finite group, and $x$ a non-Schur $p$-element of $G$. If any of the following hold:
\begin{enumerate}[\hspace{5mm} (i)]
    \item \label{item:nopsquared} $G$ contains no elements of order $p^2$.
    \item \label{item:abelianSylowCentralizer} $C_G(x)$ has an abelian Sylow $p$-subgroup.
\end{enumerate}
Then $G$ contains a extra-strong non-Schur element.
\end{lemma}
\begin{proof}
(\ref{item:nopsquared}) As $G$ has no elements of order $p^2$, neither can $C_G(x)$. Since $C_G(x)/[C_G(x),C_G(x)]$ is abelian, and has no elements of order $p^2$,  then $C_G(x)/[C_G(x),C_G(x)]$  splits as a product of a $p'$-group, and an elementary abelian $p$-group. Thus, $x \not \in E^p(C_G(x))$.

(\ref{item:abelianSylowCentralizer}) Suppose for contradiction that there are no extra-strong non-Schur elements. Take a non-Schur $p$-element $x$ of maximal order whose centralizer has an abelian Sylow $p$-subgroup, which we write $A$. Then by \cref{lem:CPx}, as $x \in E^p(C_G(x))$, we have $x \in E^p(A) = A^p$. But then $x = y^p$ for some $y\in A$, where notably $A \subseteq C_G(y) \subseteq C_G(x)$. But then $y$ has an abelian Sylow $p$-subgroup, and is of greater order than $x$, a contradiction.
\end{proof}

The following three propositions will play a crucial role in the proof of \cref{thm:g7}:

\begin{proposition}\label{prop:splitabelian}
Suppose $x$ is a $p$-element and $C_G(x)$ admits an abelian summand $A$ containing $x$. Then $G$ contains a extra-strong non-Schur element.
\end{proposition}
\begin{proof}  Take $y$ a $p$-element of maximal order so that $A$ is an abelian summand of $C_G(y)$. If $y = z^p$ for $z \in A$, then $A \subseteq C_G(z) \subseteq C_G(y)$, and $A$ is an abelian summand of $C_G(z)$, violating maximality of the order of $y$. So then $y \not \in E^p(A) = A^p$, and so clearly $y \not \in E^p(C_G(y))$. 
\end{proof}

\begin{proposition}\label{prop:centralreduction}
Suppose $H$ is a group, and $G:= H/Z(H)$. Suppose $x \in H$ is a $p$-element, $A$ is an abelian summand of $C_H(x)$ containing $x$, and the projection $C_H(x) \rightarrow C_{G}(xZ(H))$ is an isomorphism. Then $G$ contains a extra-strong non-Schur element.
\end{proposition}
\begin{proof}
    By the above proof, we know there is $y \in A$ so that $y$ is extra-strong non-Schur in $H$. Crucially this yields a group homomorphism $C_H(y) \rightarrow k^+$ on which $y$ does not vanish. This yields a homomorphism $C_H(x) \rightarrow k^+$ in the following manner:
    $$C_H(x) \rightarrow A \rightarrow C_H(y) \rightarrow k^+.$$
    Now then, utilising the inverse of the isomorphism $C_H(x) \rightarrow C_G(xZ(H))$, we obtain a homomorphism $C_G(xZ(H)) \rightarrow k^+$ not vanishing on $yZ(H)$. Restricting this to $C_G(yZ(H))$ demonstrates that $yZ(H)$ is extra-strong non-Schur.
\end{proof}
\begin{notation}

We briefly define what we mean by finite group of Lie type (since there is some ambiguity in the term \cite{MathOverflow}). Let $r$ be a prime number and let $\mathrm{GL}_n(\overline{\mathbb{F}}_r)$ be general linear group of invertible $n\times n$-matrices over $\overline{\mathbb{F}}_r$. Let ${\bf G}$ be a linear algebraic group defined over $\overline{\mathbb{F}}_r$. By  \cite[Theorem 1.7]{MalleTesterman} there is an embedding $i \colon {\bf G} \hookrightarrow \mathrm{GL}_n(\overline{\mathbb{F}}_r)$ of ${\bf G}$ as a closed subgroup into $\mathrm{GL}_n(\overline{\mathbb{F}}_r)$ for some $n$. Let $q = r^a$ for some positive integer $a$. Let $F_q \colon \mathrm{GL}_n(\overline{\mathbb{F}}_r) \to \mathrm{GL}_n(\overline{\mathbb{F}}_r)$ be the homomorphism of linear algebraic groups given by $F_q((x_{ij})) = (x_{ij}^q)$ for every  $(x_{ij}) \in \mathrm{GL}_n(\overline{\mathbb{F}}_p)$.   An endomorphism of a linear algebraic group $F \colon {\bf G} \to {\bf G}$ is called a \emph{standard Frobenius}  with respect to an $\mathbb{F}_q$-structure  if there exists a $q$ such that $i \circ F = F_q \circ i$.  We say that an endomorphism of a linear algebraic group $F \colon {\bf G} \to {\bf G}$ is a \emph{Frobenius morphism}  with respect to an $\mathbb{F}_q$-structure (also called a Steinberg endomorphism \cite[Definition 21.3]{MalleTesterman}) if there exists a positive integer $m$ such that $F^m$ is a standard Frobenius morphism. 

Now assume that ${\bf G}$ is a connected reductive algebraic group defined over $\overline{\mathbb{F}}_r$. 
Let \(F\colon \mathbf{G}\to\mathbf{G}\) be a Frobenius endomorphism of  $\mathbf{G}$. We denote the fixed-point subgroup by
\[
\mathbf{G}^{F}:=\{\,g\in \mathbf{G}\mid F(g)=g\,\},
\]
and we refer to groups of the form \(\mathbf{G}^{F}\) as \emph{finite groups of Lie type}. We say that $r$ is the \emph{defining characteristic} of $G=\mathbf{G}^{F}$. 

Let $k$ be a field. An algebraic group over $k$ is a \emph{torus} if it becomes isomorphic
to a finite product of copies of the multiplicative group $\mathbb{G}_m$ over some field containing $k$. A torus over $k$ is \emph{split} if it is isomorphic to a product of copies of $\mathbb{G}_m$ over $k$.  A reductive group is split if it contains a split maximal torus.

\end{notation}
\begin{proposition}\label{prop:liestrengtheningPair} Suppose $\mathbf{G}$ is a connected reductive algebraic group defined over $\bar{\mathbb{F}}_r$, equipped with a Frobenius endomorphism $F \colon \mathbf{G} \rightarrow \mathbf{G}$, and an embedding $\mathbf{G} \rightarrow \textrm{GL}_n(\bar{\mathbb{F}}_r)$. Suppose there $x \in \mathbf{G}$ whose embedding in $\textrm{GL}_n(\bar{\mathbb{F}}_r)$ is semisimple, and whose eigenvalues are of multiplicity one, apart from $1 \in {\bar{\mathbb{F}}_r}$. Then there is a subgroup $D \subseteq \textrm{GL}_n(\bar{\mathbb{F}}_r)$ so that $D$ is an abelian summand of $C_{\textrm{GL}_n(\bar{\mathbb{F}}_r)}(x)$ containing $x$. Therefore if $x \in \mathbf{G}^F$, then $D \cap \mathbf{G}^F$ is an abelian summand of $C_{\mathbf{G}^F}(x)$ containing $x$.
\end{proposition}

\begin{proof} Write $E_x \subseteq \bar{\mathbb{F}}_r^\times$ for the set of eigenvalues of $x$, and decompose $\bar{\mathbb{F}}_r^n = V_{E_x \backslash \{1 \}} \oplus V_1$, into eigenspaces. Then $C_{\textrm{GL}_n(\bar{\mathbb{F}}_r)}(x) = D \oplus \operatorname{GL}(V_1)$, where $D$ is the subspace of diagonal matrices on $V_{E_x\backslash{\{1\}}}$, with respect to the basis of eigenvalues of $x$. 
\end{proof}

\begin{theorem}\label{thm:HH1blockex}
    Let $k$ be a field of characteristic $p$, let $G$ be a finite group, and let $B$ be a block of $kG$ with a defect group $P$. In the following cases $\HH^1(B)\neq 0$:
    \begin{enumerate}
    \item\label{item_map_nz} The map ${\rm H}_1(P,\mathbb{F}_p)\to {\rm H}_1(G,\mathbb{F}_p)$ is nonzero.
    \item\label{item_transfer} There is $x\in P\smallsetminus E^p(P)$ such that $C_P(x)$ is a Sylow subgroup of $ C_G(x)$. 
        \item \label{item_nonSchurelt} $P$ is a Sylow subgroup of $G$, and $G$ contains a extra-strong non-Schur element.

    \end{enumerate}
In the following cases, every non-semisimple block $B$ of $kG$ satisfies $\HH^1(B)\neq 0$:
    \begin{enumerate}\setcounter{enumi}{3}
    \item\label{p_regular} For every $p$-regular element $g\in G$ such that $p$ divides $|C_G(g)|$, there is a $p$-element $x\in G$ which commutes with $g$ such that $[x]\neq 0$ in ${\rm H}_1(C_G(x),\mathbb{F}_p)$.

    \item\label{item_FGLTp} $G$ is either ${\bf G}^F$ or $[{\bf G}^F, {\bf G}^F]$ where ${\bf G}$ is a reductive group over $\overline{\mathbb{F}}_r$, split over $\mathbb{F}_r$, and where $F$ is a Frobenius map and $p=r$. 
    \item \label{item_GL} $G$ is a general linear group $\mathrm{GL}_n(\mathbb{F}_q)$, defined in any characteristic  $r$ with $q=r^n$. 
    \item\label{item_normalp} There is a normal $p$-subgroup $Q\subseteq G$ inducing a nonzero map ${\rm H}_1(Q,\mathbb{F}_p)\to {\rm H}_1(G,\mathbb{F}_p)$. 
    \end{enumerate}
\end{theorem}

\begin{proof}
(\ref{item_map_nz})  The statement follows from \cref{thm:main1}.

(\ref{item_transfer})
Set $Q:=C_P(x)$ which by assumption is a Sylow $p$-subgroup of $C_G(x)$. Since $C_P(x)\subseteq P$, we have $E^p(C_P(x))\subseteq E^p(P)$. By assumption $x\notin E^p(P)$, hence $x\notin E^p(C_P(x))$. The statement follows by combining \cref{lem:CPx} with \cref{thm:main1}.

(\ref{item_nonSchurelt}) The statement follows from \cref{thm:main1}

 (\ref{p_regular}) It is a result of Brauer \cite[10C]{Br56} that every defect group $P$ of a block $B$ of  $G$ is a Sylow $p$-subgroup of the centraliser of a $p$-regular element in $G$. Therefore, we may assume that $P$ is a Sylow subgroup of $C_G(g)$, and since $x$ is a $p$-element that commutes with $g$, we may assume that $x\in P$. Then the statement follows from Theorem \ref{thm:main1}

(\ref{item_FGLTp})
If $G$ is either ${\bf G}^F$ or $[{\bf G}^F, {\bf G}^F]$ where ${\bf G}$ be a reductive group over $\overline{\mathbb{F}}_r$, split over $\mathbb{F}_r$, whose root system is irreducible, then blocks of  ${\bf G}^F$ and of $[{\bf G}^F, {\bf G}^F]$ are all either of maximal defect or of defect zero, see the main result in \cite{Humphreys71}. 
It follows from \cite[Lemma 3.1 (1)]{FJL}
that ${\bf G}^F$ satisfies $S(p)$.  By the same lemma,   ${\bf G}^F$ cointains  non-Schur $p$-elements  $x$ such that $C_G(x)$ is abelian.  By  \cref{lemma:StrongSuperStrong} \eqref{item:abelianSylowCentralizer} we have ${\bf G}^F$ also satisfies $ES(p)$ and the statement follows from (\ref{item_nonSchurelt}).

To finish, note that we may assume the root system of ${\bf G}$ is irreducible, since otherwise then ${\bf G}$ would be product of reductive groups. In this case ${\bf G}^F$ (or $[{\bf G}^F, {\bf G}^F]$) would correspondingly split as a product of groups. By considering the group algebras and the corresponding blocks, this reduces to the statement just established.

(\ref{item_GL}) Since  $\mathrm{GL}_n(k)$ is a split reductive group over any field $k$  \cite[Example 21.6]{Milne17}, by (\ref{item_FGLTp}) it is enough to consider the non-defining characteristic case. Following the proof in  \cite[Theorem 7.18]{SL2} each block of a general linear group 
is derived equivalent to the principal block of a product of general linear groups. From \cite[Proposition 3.2 (a)]{FJL}, the element they call $f_\alpha$ satisfies the criteria of \cref{prop:liestrengtheningPair}. Thus we obtain an abelian summand $A$ of $C_G(x)$ containing $x$, and \cref{prop:splitabelian} completes the claim.

(\ref{item_normalp}) Every defect group $P$ must contain $Q$ by \cite[Theorem 6.2.6]{Linckelmann18-II}, and therefore if ${\rm H}_1(Q,\mathbb{F}_p)\to {\rm H}_1(G,\mathbb{F}_p)$ is nonzero then ${\rm H}_1(P,\mathbb{F}_p)\to {\rm H}_1(G,\mathbb{F}_p)$ is nonzero. Hence this follows from (\ref{item_map_nz}). 
\end{proof}

    We end this section with a remark on \cref{thm:HH1blockex} and its precedents.

\begin{remark}\label{rmk:knowncases}
    Let $B$ be a non-semisimple block of $kG$.  Todea has shown that that $\HH^1(B)\neq 0$ if $G$ is $p$-solvable, and also that $\HH^1(B)\neq 0$ if the defect group of $B$ is normal in $G$ \cite[Corollary 1.4]{Todea}.  We point out that under the additional assumption that $B$ has a Sylow defect group, these results follow from \cref{thm:HH1blockex}. Todea's arguments rely on the Morita equivalences constructed by K\"ulshammer \cite[Theorem A]{Ku}, while (in the special case of Sylow subgroups) our arguments are more direct.

    Brou\'{e}'s abelian defect group conjecture predicts that if $B$ has an a abelian defect group $P$, then $B$ is derived equivalent to a block of $kN_G(P)$, also having $P$ as a defect group. Since $P$ is normal in $N_G(P)$, it follows from \cite[Corollary 1.4]{Todea} that, as long as $P$ is nontrivial, the conjecture predicts $\HH^1(B)\neq 0$. The conjecture is known to hold in various cases; see for example \cite[Chapter 14]{AHK} and references therein. 
    We point out that if $B$ is a block with a nontrivial abelian \emph{Sylow} defect group $P$, then the nonvanishing of $\HH^1(B)$ follows from \cref{thm:HH1blockex},  
    noting that $G$ contains a non-Schur $p$-element since $P$ is abelian.
\end{remark}

\section{Blocks of the symmetric and alternating groups}

In this section we show the validity of Linckelmann's Conjecture for blocks of the symmetric group $S_n$ and of the alternating group $A_n$. 

The case of $S_n$ this has been established already in \cite[Theorem 1.1]{BKL23} by an explicit calculation of the dimension of $\HH^1(B)$. 
 A crucial ingredient in their proof is the reduction to principal blocks of $kS_n$. Indeed, each block $B$ of the symmetric group algebra $kS_n$ is associated with a non-negative integer weight $w$, and by \cite[Theorem 7.2]{SL2}, any two blocks (possibly of different symmetric groups) with the same weight $w$ are derived equivalent. Since $\HH^1(B)$ is a derived invariant, it suffices to consider one block for each weight $w$, for instance the principal block of $kS_{pw}$.

We give a proof in the $S_n$ case that is more elementary and avoids relying on \cite{SL2}. In addition, it provides the key idea for the $A_n$ case where there is no analogue of the results in \cite{SL2}.

\begin{theorem}\label{thm_Sncase}
    Let $k$ be a field of characteristic $p$, let $G=S_n$, and let $B$ be a block of $kG$ with a nontrivial defect group $P$. Then $\HH^1(B)\neq 0$.
\end{theorem}

\begin{proof}
We will apply \cref{thm:HH1blockex} part  (\ref{p_regular}). Let $g$ be a $p$-regular element in $S_n$ with cycle type $(1)^{r_1}(2)^{r_2}\cdots (n)^{r_n}$, where $n_i$ is the number of cycles of length $i$ in the disjoint cycle representation of $g$. Since $g$ is $p$-regular, the cycles $(i)$ appearing with $r_i>0$ have length coprime to $p$.

The centraliser of $g$ is $C_G(g)= S_{r_1}\times \dots \times (C_n\wr S_{r_n})$, and we assume this has order divisible by $p$. Therefore there is an $i$, coprime to $p$, such that the factor $C_i\wr S_{r_i}$ has order divisible by $p$. This means $p$ divides $|S_{r_i}|$, so we may take  a $p$-cycle $(1 \cdots p)$ in $S_{r_i}$. 

The  $p$-element $x=(1,\dots, 1,(1\cdots p))\in  C_i\wr S_{r_i}$ is sent under the inclusion $C_i\wr S_{r_i}\hookrightarrow  S_n$ to a product of $i$ disjoint $p$-cycles. 

Therefore, the cycle type of  $x$ in $S_n$ is $(1)^{n-ip}(p)^i$. So $C_G(x)=S_{n-ip}\times C_p\wr S_i$. Consequently ${\rm H}_1(C_G(x),\mathbb{F}_p)={\rm H}_1(S_{n-ip}\times C_p\wr S_i,\mathbb{F}_p)\cong {\rm H}_1(S_{n-ip} ,\mathbb{F}_p)\times C_p \times {\rm H}_1(S_i,\mathbb{F}_p)$, in which the $p$-element $x$ is sent to $[x]=h^i\neq 1$ in $C_p=\langle h\rangle$, using the fact that $p$ does not divide $i$. The statement now follows from \cref{thm:HH1blockex} part (\ref{p_regular}).
\end{proof}

Before proceeding with the $A_n$ case we recall that to each block $B$ of $kA_n$ there is associated a non-negative integer  $w \leq \lfloor n/ p\rfloor$, called the \emph{weight} of $B$, such that the Sylow
$p$-subgroups of $A_{pw}$ are defect groups of $B$, where $A_{pw} \subseteq A_n$ is the natural inclusion. 

Given a block $B$ of $A_n$ with weight $w$ and defect group $P\subseteq A_{pw} \subseteq A_n$, there is a block $\tilde{B}$ of $S_n$ with defect group $\tilde{P}\subseteq S_{pw} \subseteq S_n$ such that  $P=\tilde{P}\cap A_n$.  See \cite{Olsson} for more information.

\begin{theorem}\label{thm_Ancase}
    Let $k$ be a field of characteristic $p$, let $G=A_n$ be the alternating group, and let $B$ be a block of $kG$ with a nontrivial defect group $P$. Then $\HH^1(B)\neq 0$.
\end{theorem}

\begin{proof} We may write $P = \tilde P\cap A_n$ where $\tilde P$ is the defect group of a block $\widetilde{B}$ of $S_n$. In the proof of \cref{thm_Sncase} we constructed $x \in \tilde P$ so that $x \not \in E^p(C_{S_n}(x))$. Notably, for odd $p$, this element lies in $A_n$. Thus $x \in P$, and as $[C_{S_n}(x):C_{A_n}(x)]$ is coprime to $p$, by \cref{lem:CPx} we know $x$ is an extra-strong non-Schur element of $A_n$.

For $p=2$ the  proof makes use of weight theory. We divide in two cases:  $w> 3$ and  $0<w<3$, where $w$ is the weight of a block $B$ of $A_n$.

    For $w\geq 3$, since the defect group $P$ of $B$ is isomorphic to a Sylow $2$-subgroup of $A_{2w}$, we can consider an element $x\in P$ which has cycle structure $(4,2,1^{n-6})$, say $x=(1234)(56)$. Then 
    $C_{S_n}(x)= C_4 \times C_2 \times S_{n-6}$ and $[x]$ is non-zero in $\mathrm{H}_1(C_{S_n}(x),k)$. Looking at the map $\mathrm{H}_1(C_{A_n}(x), k) \to \mathrm{H}_1(C_{S_n}(x), k)$, it must be that $[x]$ is non-zero in $\mathrm{H}_1(C_{A_n}(x),k)$, as desired.

      If $0<w<3$, we need to only consider  $w=2$ since if $w=1$, then $A_2$ is the trivial group. If $w=2$, then $\tilde{P}$ is a Sylow $2$-subgroup of $S_4$. Therefore $B$ has as a defect group $P\cong C_2\times C_2$. By \cite[Corollary 2]{Linckelmann94}  $B$ has either one or three simple modules. In the former case $B$ is Morita equivalent to $kP$. In the latter case, $B$ is derived equivalent to $kA_4$. The statement follows.
\end{proof}

\begin{remark}
In the proof of \cref{thm_Ancase} for the  $w=2$ case, we cannot pick an element $x\in P \cong C_2 \times C_2$ such that $[x]\neq {\rm H}_1(C_{A_n}(x), \mathbb{F}_2)$ for an arbitrary $n$.  Indeed, if $n\geq 6$, it is easy to check that for any $x\in P$ we have $x\in [C_{A_n}(x), C_{A_n}(x)]$.  This shows that \cref{thm:main1} provides a sufficient but not necessary condition for the non-vanishing of $\HH^1(B)$.
\end{remark}

\section{The extra-strong non-Schur property for primes greater than five}\label{sec_char_5}

   Fleischmann, Janiszczak, and Lempken prove that every finite group satisfies the weak non-Schur property $W(p)$ for every prime \cite{FJL}. This was established by reducing to  the case of finite simple groups and using the classification. As a consequence of the centraliser decomposition for the Hochschild cohomology of group algebras, one obtains from this that $\HH^1(kG) \neq 0$ for every finite group with order divisible by $p={\rm char}(k)$. 

We saw in \cref{thm:HH1blockex} that if the extra-strong non-Schur property holds, then blocks $B$ of $kG$ of maximal defect satisfy $\HH^1(B)\neq 0$. We have also seen that unfortunately the strong non-Schur property $S(p)$ does not hold for every group and every prime; by observing that  $S(p)$ implies $SC(p)$ and using  \cref{rmk:Scp}. So in particular the extra-strong non-Schur property $S(p)$ also does not hold. However, using similar methods to \cite{FJL} we are able to show that $ES(p)$ holds whenever $p$ is bigger than $5$.

\begin{theorem}\label{thm:g7}
If $p$ is a prime number greater than $5$ then any finite group $G$ with order divisible by $p$ contains a extra-strong non-Schur $p$-element. In particular, if $k$ is a field of characteristic $p$, then any block $B$ of $kG$ of maximal defect satisfies $\HH^1(B)\neq 0$.
\end{theorem}

Before starting the the proof we will fix some notation and recall a useful result on simple Lie--Chevalley groups.

\begin{notation}
Let ${\bf G}$ be a connected reductive algebraic group defined over $\overline{\mathbb{F}}_r$, and let ${\bf T}\subseteq {\bf G}$ be a maximal torus. The group
\[
W:=N_{\bf G}({\bf T})/C_{\bf G}({\bf T})
\]
is called the Weyl group  of $\bf{ G}$ with respect to ${\bf T}$.  

A \emph{finite simple Lie--Chevalley group} $G$ is a finite simple group which is either a finite group of Lie type ${\bf G}^F$, or a commutator subgroup or central factor of some ${\bf G}^F$. In either case we say  that the group $W$ above is the \emph{Weyl group} of $G$ (leaving ${\bf G}$ implicit). 
\end{notation}

\begin{lemma}\cite[Lemma 3.1 (2)]{FJL}\label{lemma:FJL2}
Let $G$ be a finite simple Lie--Chevalley group defined over a field of characteristic $r$,  
with Weyl group $W$. If 
$p\nmid |W|$ is a prime different from $r$, then the Sylow $p$-subgroups of $G$ are abelian. \qed
\end{lemma}

\begin{proof}[Proof of \cref{thm:g7}.]
By \cref{thm:HH1blockex} (\ref{item_nonSchurelt}) and   \cref{Prop:red} it suffices to prove that every finite simple group satisfies $ES(p)$ for $p>5$. More than this, we will make an effort to decide for each finite simple group those primes (including $p=2,3,5$) for which $ES(p)$ holds (although we stop short of a complete answer on this problem). In many cases, we  show that a finite simple group  satisfies $S(p)$, by relying on  computation essentially contained in \cite{FJL}, and then we will apply Lemma \ref{lemma:StrongSuperStrong} or  
\cref{prop:splitabelian}. 
We will first list the cases already covered by \cite{FJL}, and then we finish the proof by analysing the remaining cases directly.

{\bf Sporadic simple groups.} By \cite[Proposition 2.2]{FJL}, if $G$ is one of the $26$ sporadic simple groups then $G$ satisfies $S(p)$ with the exception of $(G,p)$ being $({\rm Ru}, 3)$, $({\rm J}_4,3)$ or $({\rm Th},5)$ where ${\rm Ru}$, ${\rm J}_4$, and ${\rm Th}$ denote the Rudvalis group, the Janko group ${\rm J}_4$, and the Thompson group, respectively. For this reason we exclude $p=3$ and $5$. For $p>5$ there are no elements of order $p^2$, then the  statement follows from   \cref{lemma:StrongSuperStrong} \eqref{item:nopsquared}.

{\bf Cyclic groups.} It is straightforward to see that the cyclic groups satisfy $ES(p)$ for all $p$.

{\bf Alternating groups.}  In  the proof of \cref{thm_Ancase} it has been shown that $A_n$  satisfies  $ES(p)$ for $p>2$.

{\bf Simple Lie--Chevalley groups  in defining characteristic.} By \cite[Lemma 3.1]{FJL} $G$  contains a regular unipotent element, and any such element $x$ has order a power of $p$ and has abelian centraliser $C_G(x)$; in particular $x$  is a non-Schur $p$-element. By  \cref{lemma:StrongSuperStrong} \eqref{item:abelianSylowCentralizer} $G$ also satisfies $ES(p)$.

{\bf Simple Lie--Chevalley groups in non-defining characteristic.} For the classical groups we follow the notation in  \cite[Chapter 2]{Atlas}. For the rest of the simple  groups of Lie type we follow \cite[Chapter 3]{Atlas}.  Let $p,r$ be two different primes and let $q$ be a power of  $r$.

Firstly, let $G$ be one of the following Lie--Chevalley groups of classical type:
    \begin{multicols}{3}
    \begin{enumerate}   
    \item $A_n(q)=\mathrm{PSL}_n(q)$ 
    \item $^2A_n(q^2)=\mathrm{U}_{n+1}(q)$ 
    \item  $B_n(q)=\mathrm{O}_{2n+1}(q)$ 
    \item $C_n(q)=\mathrm{PSp}_{2n}(q)$ 
    \item $D_n(q)=\mathrm{O}_{2n}^+(q)$
    \item $^2D_n(q^2)=\mathrm{O}_{2n}^{-}(q)$ 
 \end{enumerate}
 \end{multicols}
for the appropriate values of $n$ and $q$  to make $H$ simple. Then $H$ satisfies $S(p)$ for all $p\neq r$ by \cite[Proposition 3.2]{FJL}. 
We seek to obtain a extra-strong non-Schur element from their proof. To do this, it is worth highlighting their strategy of proof. For each such group $G$, we fix a group $H \subseteq \operatorname{GL}_n(\bar{\mathbb{F}}_r)$ for which $G = H/Z(H)$ or $G =H$, in the standard manner (in the cases above these are $\textrm{SL}_n(q), \textrm{SU}_n(q), \textrm{SO}_n(q), \textrm{Sp}_n(q), \Omega_n^\epsilon(q)$). Then, an element $x \in H$ (which they write as $f_\alpha$) is chosen obeying the hypotheses of \cref{prop:liestrengtheningPair}. They finish by demonstrating that $C_H(x) \rightarrow C_G(xZ(H))$ is an isomorphism in all relevant cases. So combining \cref{prop:liestrengtheningPair} with \cref{prop:centralreduction} we are done. 
We next consider the following simple Lie--Chevalley groups, for the appropriate values of $n$ and $q$:  
    \begin{multicols}{3}
    \begin{enumerate} 
    \setcounter{enumi}{6}
        \item $^2B_2(2^{2n+1})$ 
    \item $^2G_2(3^{2n+1})$ 
    \item \label{item_2F4} $^2F_4(2^{2n+1})$ 
    \item \label{item_G2}$G_2(q)$ 
    \item \label{item_F4} $F_4(q)$
    \item \label{item_3D}$^3D_4(q^3)$
    \item \label{item_E6}$E_6(q)$
    \item \label{item_E62}$^2E_6(q^2)$
 \end{enumerate}
 \end{multicols}
   \begin{enumerate}[$\bullet$]
   \item The order of the Weyl group of  $^2B_2(2^{2n+1})$ is $2^3$ \cite[Table 2]{Hump90}. As we are in non-defining characteristic  $p\neq r=2$, \cref{lemma:FJL2} together with \cref{prop_nonSchur}\eqref{item_sylabRel} imply that $ES(p)$ holds for $^2B_2(2^{2n+1})$. 

\item For $^2G_2(3^{2n+1})$ the order of the Weyl group is $2^2\cdot3$ \cite[Table 2]{Hump90}, and $p\neq r=3$, so by Lemma \ref{lemma:FJL2} we only need to consider $p=2$. All Sylow $2$-subgroups  of  $^2G_2(3^{2n+1})$ are abelian \cite[Proof of Proposition 4.1]{FJL} therefore   $^2G_2(3^{2n+1})$  satisfies $ES(p)$ for all $p$ by  \cref{prop_nonSchur}\eqref{item_sylabRel}.
   \item For $^2F_4(2^{2n+1})$  we need to exclude $p=3$ since the order of the Weyl group is $2^7\cdot 3^2$ \cite[Table 2]{Hump90} and $p\neq r=2$.
       \item For $G_2(q)$, $F_4(q)$ and $^3D_4(q^3)$ we need to exclude $p=2,3$ as the order of the Weyl groups are $2^2\cdot 3$, $2^7\cdot 3^2$and $2^6\cdot 3$ \cite[Table 2]{Hump90}, respectively.
       \item For $E_6(q)$, $^2E_6(q^2)$, we need to exclude $p=2,3,5$ since the order of the Weyl group is $2^7\cdot 3^4\cdot 5$ \cite[Table 2]{Hump90}.   
   \end{enumerate}   
Only three cases remain, and for these we will need to do our own analysis:
    \begin{multicols}{3}
    \begin{enumerate} 
    \setcounter{enumi}{14}
    \item \label{item_Tits}$^2F_4(2)'$ (Tits group)
    \item \label{item_E7}$E_7(q)$
    \item \label{item_E8}$E_8(q)$
 \end{enumerate}
 \end{multicols}
\begin{enumerate}[$\bullet$]
\item  The Tits group $^2F_4(2)'=[{}^2F_4(2),{}^2F_4(2)]$ has order $2^{11}\cdot 3^3\cdot 5^2 \cdot 13$ \cite[p74]{Atlas}. It has been observed in \cite[p3]{FJL} that $G$ does not satisfy $S(3)$, so in particular it does not satisfy $ES(3)$.
     For $p=5$ and $13$, the Sylow subgroups are abelian since their order is at most $p^2$. For $p=2$ consider an element $x$ of order $16$ in the conjugacy class $16A$ (using the ATLAS notation \cite{AtlasV3}).  Since the group generated by $x$ has order $16$ and is contained in $C_G(x)$, we have that $C_G(x)\cong C_{16}$ is abelian hence $G$ satisfies $S(2)$, and so $ES(2)$ holds by \cref{lemma:StrongSuperStrong}. In summary, $G$ satisfies $ES(p)$ for all primes excluding $p=3$.
     
\item  Let $G=E_7(q)$. Since the order of the Weyl group is $2^{10}\cdot 3^4\cdot 5\cdot 7$ \cite[Table 2]{Hump90}, $G$ satisfies $S(p)$ when $p\neq 2,3,5,7$. We have already abandoned hope on $p=2,3,5$, so it remains to show that $ES(7)$ holds. Following the proof of \cite[Proposition 2.5]{NT}, for $p=7$  the Sylow subgroups of $G$ are isomorphic to $P_1=C_{7^c}\wr  C_7$.  We have $[P_1,P_1]\subseteq (C_{7^c})^{7}$, as $(C_{7^c}\wr  C_7)/ (C_{7^c})^{7}\cong C_7$ is abelian. Consequently $\mathrm{exp}([P,P])\leq 7^c$. Let $x=((g,1,\ldots,1), h)\in C_{7^c}\wr  C_7$, where $g, h$ is are generators of $C_{7^c}$ and $C_7$, respectively. Then $x^7=(g,\dots,g, 1)$  has order $7^c$. Hence $\mathrm{exp}(P)\geq 7^{c+1}$. The statement follows from \cref{prop_nonSchur}\eqref{item_syl_expRel}. 

 \item Let $G=E_8(q)$. Since the order of the Weyl group is $2^{14}\cdot 3^5\cdot 5^2\cdot 7$ \cite[Table 2]{Hump90}, $G$ satisfies $S(p)$ when $p\neq 2,3,5,7$. This time, for $p=7$, the proof of \cite[Proposition 2.5]{NT} shows that  the Sylow subgroups of $G$ are isomorphic to  $P=P_1\times P_2$ where $P_1=C_{7^c}\wr  C_7$ and $P_2=C_{7^c}$. So by \cref{lem:CPx}  it is enough to take a generator  $x\in P_2$ and observe $(1,x)\in Z(P)$ but $(1,x)\notin E^p(P,P)$. 
\end{enumerate}
This finishes the proof that every finite simple group satisfies $S(p)$ for $p>5$.
\end{proof}

As a consequence of the proof of \cref{thm:g7} we obtain an interesting result in group theory.

\begin{corollary}
    Let $p$ be a prime number greater than $5$. Then  $S(p)$ implies $ES(p)$.
\end{corollary}

\begin{proposition}
    All sporadic simple groups, apart from $(\mathrm{Ru},3)$, $(\mathrm{Th},5)$ and $(J_4,3)$,  have the extra-strong non-Schur property $ES(p)$.
\end{proposition}
\begin{proof}
      Let $G$ be sporadic simple group  and $p$ a prime number. Consider the pair $(G,p)$  apart from $(\mathrm{Ru},3)$, $(\mathrm{Th},5)$ and $(J_4,3)$. By the proof of \cref{thm:g7} we only need  to consider $p\leq 5$. For all groups excluding $\mathrm{B}, \mathrm{M}$ and $\mathrm{Th}$, direct computation in GAP produces extra non-Schur elements. In the remaining cases we proceed as follows: we search for a $p$-element of maximal order in $G$. Then use the ATLAS to identify the order of its centralizer, and search all groups with this property using GAP. Finally, we verify that all central elements of the correct order are not contained in $E^p(G)$. For example, for $(Th,2)$ the maximal order of a $2$-element is $8$, and its centralizer is of order $96$. There are $8$ groups of order $96$ and exponent $8$ admitting a central element of order $8$, and in all cases, the central elements of order $8$ are not contained in the $2$-Frattini.
\end{proof}
\begin{remark}
    In a forthcoming paper the first and the third authors, in joint work with Can Wen, investigate the dimensions of \(\HH^{1}(B)\) for blocks of sporadic finite simple groups.
\end{remark}

\section{Triangular algebras and separable equivalences}
\label{sec_Hochschild}

In this section $k$ is an arbitrary field.

\begin{definition}
Let $A$ and $B$ be algebras and let $M$ be an $A$-$B$-bimodule, that is, a $A^{\rm op}\otimes B$-module. The triangular algebra of $M$ is the matrix algebra
\[
T_M\colonequals\begin{pmatrix}
    A & M\\
    0 & B
\end{pmatrix}.
\]
In other words, $T_A=A\oplus M\oplus B$ with its algebra structure inherited from that of $A$ and $B$ and from the bimodule structure of $M$, and with $M$ forming a two-sided, square-zero ideal.
\end{definition}

Inside the triangular algebra there are two idempotents
\[
e_A\colonequals\begin{pmatrix}
    1 & 0\\
    0 & 0
\end{pmatrix}\quad\text{and}\quad e_B\colonequals\begin{pmatrix}
    0 & 0\\
    0 & 1
\end{pmatrix},\]
and we will use these idempotents to define restriction homomorphisms
\[
\HH^*(A,A)\xleftarrow{\rho^M_A}\HH^*(T_M,T_M)\xrightarrow{\rho^M_B} \HH^*(B,B).
\]

\begin{definition}\label{def_restriction}
Note that if $P\to T_M$ is a $T_M$-$T_M$-bimodule resolution of $T_M$ then $e_APe_A\to e_ATe_A=A$ is an $A$-$A$-bimodule resolution of $A$. Because of this, we may represent a class $[f]\in \HH^n(T_M,T_M)$ as a chain map $f\colon P\to T_M[n]$ and define
\[
\rho^M_A(f)\colonequals    e_Afe_A\colon e_APe_A\longrightarrow e_AT_Me_A[n]=A[n] \quad \text{in }\HH^n(A,A).
\]
Likewise, $\rho^M_B$ is defined using the idempotent $e_B$. These restriction homomorphisms were used by Keller in \cite{DIH} 
to establish invariance properties of the higher structure on Hochschild cohomology. 
\end{definition}

We will also make use of tensor homomorphisms 
\[
\HH^*(A,A)\xrightarrow{\tau^M_A}\Ext^*_{A^{\rm op}\otimes B}(M,M)\xleftarrow{\tau^M_B} \HH^*(B,B).
\]
\begin{definition}
Assume that $M$ is, separately, projective as a left $A$-module and a right $B$-module. Note that if $Q\to A$ is a $A$-$A$-bimodule resolution of $A$ then $Q\otimes_AM \to A\otimes_AM=M$ is an $A$-$B$-bimodule resolution of $M$. Because of this, we may represent a class $[f]\in \HH^n(A,A)$ as a chain map $f\colon Q\to A[n]$ and define
\[
\tau^M_A(f)\colonequals    f\otimes_AM \colon Q\otimes_A M \longrightarrow A\otimes_A M [n]=M[n] \quad \text{in }\Ext^*_{A^{\rm op}\otimes B}(M,M).
\]
Likewise, $\tau^M_B$ is defined by applying the exact functor $M\otimes_B-$. These tensor homorphisms are often used to define transfer maps on Hochschild cohomology, c.f.~\cite{transfer}.
\end{definition}

One of our main tools is the following Mayer--Vietoris exact sequence. 

\begin{theorem}[Keller]\label{th_MV}
Let $M$ be an $A$-$B$-bimodule that is projective as a left $A$-module and a right $B$-module. There is an exact sequence of restriction and tensor homomorphisms
\[
 \cdots\to \HH^*(T_M,T_M)\to \HH^*(A,A)\oplus \HH^*(B,B) \to \Ext^*_{A^{\rm op}\otimes B}(M,M)\to \HH^{*+1}(T_M,T_M)\to \cdots
\]
In particular $\im(\rho^M_A)=(\tau^M_A)^{-1}({\rm im}(\tau^M_B))$ and $\im(\rho^M_B)=(\tau^M_B)^{-1}({\rm im}(\tau^M_A))$.
\end{theorem}

\begin{proof}
    The existence of the long exact sequence is  \cite[(4.5.2)]{DIH}, and the stated equalities follow by diagram chasing.
\end{proof}

\begin{remark}\label{rem_bimod}
    If $\varphi\colon A\to B$ is a algebra homomorphism making $M=B$ into an $A$-$B$-bimodule, and if we assume that $B$ is projective as a left $A$-module, then by the derived tensor-hom adjunction
    \[
    \Ext^*_{A^{\rm op}\otimes B}(B,B)\cong \HH^*(A,B).
    \]
    Making this identification, $\tau^B_B$ corresponds to the restriction of scalars $\HH^*(B,B) \to \HH^*(A,B)$ along $\varphi$; in degree one this is modeled on derivations as the map $\Der(B,B)\to \Der(A,B)$ sending $\theta$ to $\theta\varphi$.
    
    Likewise, $\tau^B_A$ corresponds to the extension of coefficients $\HH^*(A,A) \to \HH^*(A,B)$ along $\varphi$, thought of as an $A$-$A$-bimodule homomorphism; in degree one this is modeled on derivations as the map $\Der(A,A)\to \Der(A,B)$ sending $\theta$ to $\varphi\theta$.
    
\end{remark}

\begin{definition}\cite{Kadison}
Let $M$ be an $A$-$B$-bimodule and $N$ be a $B$-$A$-bimodule. We say that $(M,N)$ induces a separable equivalence 
if $A$ is a direct summand of $M\otimes_B N$ as an $A$-$A$-bimodule, and $B$ is a direct summand of $N\otimes_A M$ as an $B$-$B$-bimodule.
\end{definition}

\begin{theorem}\label{th_sep}
If $(M,N)$ induces a separable equivalence then $\tau^M_A$ and $\tau^M_B$ are injective, and the two restriction homomorphisms
\[
\HH^*(A,A)\xleftarrow{\rho^M_A}\HH^*(T_M,T_M)\xrightarrow{\rho^M_B} \HH^*(B,B)
\]
have isomorphic images.
\end{theorem}

\begin{proof}
    Suppose that $\iota\colon A\to M\otimes_BN$ and $\pi \colon M\otimes_BN\to A$ are $A$-bilinear maps satisfying $\pi\iota=1$. Then the composition
    \[
    \HH^*(A,A) \xrightarrow{\tau^M_A} \Ext^*_{A^{\rm op}\otimes B}(M,M) \xrightarrow{-\otimes_BN} \Ext^*_{A^{\rm op}\otimes A}(M\otimes_BN,M\otimes_BN) \xrightarrow{\pi(-)\iota} \Ext^*_{A^{\rm op}\otimes A}(A,A)
    \]
    is the identity on $\HH^*(A,A)$, and so $\tau^M_A$ is injective. Likewise, $\tau^M_B$ is injective since $B$ is a summand of $N\otimes_AM$.

    For the second statement, it is equivalent to prove that the two homomorphisms have the same kernel. So we will show, for $x\in \HH^*(T_M,T_M)$, that $\rho_A^M(x)=0$ implies $\rho_B^M(x)=0$, the converse implication following by symmetry. Using the Mayer--Vietoris sequence above $\tau^M_B\rho_B^M(x)=\tau^M_A\rho_A^M(x)=0$, and so $\rho_B^M(x)=0$ since $\tau^B_M$ is injective.
\end{proof}

For the rest of the section we specialise to the case of a group algebra $kG$, and we take $B$ to be a block of $kG$ with defect group $P$. In particular we will consider the triangular algebras
\[
T_G=\begin{pmatrix}
    kP & kG\\
    0 & kG
\end{pmatrix}\quad\text{and}\quad
T_B=\begin{pmatrix}
    kP & B\\
    0 & B
\end{pmatrix},
\]
viewing $kG$ as a $kP$-$kG$-bimodule, and $B$ as a $kP$-$B$-bimodule in the natural way.

\begin{lemma}\label{lem_fac}
If $B$ is a block of $kG$ with defect group $P$ then  
$\rho^{kG}_{kP}$ factors as
\[
\HH^*(T_G,T_G)\xrightarrow{\ \ \ }\HH^*(T_B,T_B)\xrightarrow{\rho^B_{kP}} \HH^*(kP,kP).
\]
\end{lemma}

\begin{proof}
The block idempotent $b\in kG$ with $B=b(kG)$ yields an idempotent $e_b=\begin{pmatrix}
    1&0\\0&b
\end{pmatrix}$ in $ T_G$ such that $T_B=e_bT_Ge_b$. Just as in \cref{def_restriction} we may define 
\[
\HH^*(T_G,T_G)\to\HH^*(T_B,T_B),\quad f\mapsto e_bfe_b,
\]
and the claimed factorisation holds because $\begin{pmatrix}
    1&0\\0&b
\end{pmatrix} \begin{pmatrix}
    1&0\\0&0
\end{pmatrix} =\begin{pmatrix}
    1&0\\0&0
\end{pmatrix} \begin{pmatrix}
    1&0\\0&b
\end{pmatrix} = \begin{pmatrix}
    1&0\\0&0
\end{pmatrix}$.
\end{proof}

Combining the above ingredients yields the next result, which includes the statement claimed in the introduction,  just before \cref{thm_1_cor}, as the degree one case.

\begin{theorem}\label{thm_A_general}
For any block $B$ of $kG$ with defect group $P$, using the natural homomorphisms
\[
\HH^*(kP,kP)\xrightarrow{\alpha}\HH^*(kP,kG)\xleftarrow{\beta} \HH^*(kG,kG)
\]
there is an embedding  $\alpha^{-1}({\rm im}(\beta))\hookrightarrow\HH^*(B,B)$.
\end{theorem}

\begin{proof}
 The pair $(B,B)$, considered as a $B$-$kP$- and a $kP$-$B$-bimodule, induces a separable equivalence by the definition (according to taste) of defect groups \cite[3.1]{LinckelmannEMS}. Therefore $\im(\rho_{kP}^B)\cong \im(\rho_{B}^B)\subseteq \HH^*(B,B)$ by \cref{th_sep}. Moreover $\im(\rho_{kP}^B)\supseteq \im(\rho_{kP}^{kG})$ by \cref{lem_fac}, and  $\im(\rho^{kG}_{kP})=(\tau^{kG}_{kP})^{-1}({\rm im}(\tau^{kG}_{kG}))$ by \cref{th_MV}. Finally, by \cref{rem_bimod}, we may identify $\alpha=\tau^{kG}_{kP}$ and $\beta=\tau^{kG}_{kG}$.
\end{proof}

\begin{proposition}\label{prop_centraliser_decomp}
Let $H$ be a subgroup of $G$, and consider the natural homomorphisms
\[
\HH^1(kH,kH)\xrightarrow{\alpha}\HH^1(kH,kG)\xleftarrow{\beta} \HH^1(kG,kG).
\]
Then $\alpha^{-1}({\rm im}(\beta))\neq 0$ if and only if there is an element $x\in H$ such that the homomorphism
\[
{\rm H}_1(C_H(x),\mathbb{F}_p)\to {\rm H}_1(C_G(x),\mathbb{F}_p)
\]
induced by the inclusion $C_H(x)\subseteq C_G(x)$ is nonzero and for any element $y=hxh^{-1}$ such that $y\notin H$  the homomorphism
\[
{\rm H}_1(C_H(y),\mathbb{F}_p)\to {\rm H}_1(C_G(y),\mathbb{F}_p)
\]
is zero.
\end{proposition}

\begin{proof}
We use the centraliser decomposition of Hochschild cohomology \cite[Theorem 2.11.2]{Benson_1991}. 

A slightly extended version of the usual statement (with the same proof) yields a commuting diagram
\[
\begin{tikzcd}
\HH^1(kH,kH)\ar[r,"{\alpha}"] \ar[d,"\cong"] & \HH^1(kH,kG) \ar[d,"\cong"]  & \HH^1(kG,kG) \ar[l,"{\beta}"'] \ar[d,"\cong"] \\
\bigoplus_{[x]\in H/{{\sim
}_H}} {\rm H}^1(C_H(x),k)\ar[r] & \bigoplus_{[x]\in G/\sim_H} {\rm H}^1(C_H(x),k)& \bigoplus_{[x]\in G/{{\sim
}_G}} {\rm H}^1(C_G(x),k). \ar[l]
\end{tikzcd}
\]
Here, the left-hand sum is indexed by a set of representatives in $H$ for the set $H/{\sim_H}$ of conjugacy classes in $H$, and similarly for the right-hand sum, while the middle sum is indexed by representatives for the set $G/\sim_H$ of equivalence classes for the conjugation action of $H$ on $G$. The lower left map is the inclusion of those summands indexed by the subset $H/{\sim_H}\subseteq G/{\sim_H}$. To describe the lower right map, any class $[x]\in G/{\sim_G}$ is a disjoint union $[x]=[x_1]\cup \cdots \cup[x_n]$ of classes in $G/\sim_H$, and each $x_i=g_ixg_i^{-1}$ for some $g_i\in G$; the map is then given by the sum over $i=1,\ldots,n$ of the restriction homomorphisms 
\[
{\rm H}^1(C_G(x),k) \xrightarrow{(g_i^{-1}(-)g_i)^*} {\rm H}^1(C_G(x_i),k)\xrightarrow{{\rm res}^{C_G(x_i)}_{C_H(x_i)}} {\rm H}^1(C_H(x_i),k).
\]
From this description it follows that $\alpha^{-1}({\rm im}(\beta))\neq 0$ if and only if there is a summand of the lower row
\[
{\rm H}^1(C_H(x),k)\xrightarrow{\ \ =\ \ } {\rm H}^1(C_H(x),k) \xleftarrow{{\rm res}^{C_G(x)}_{C_H(x)}} {\rm H}^1(C_G(x),k),
\]
indexed by some $x\in H$, with ${\rm res}^{C_G(x)}_{C_H(x)}\neq 0$ and for any element $y=hxh^{-1}$ such that $y\notin H$  the homomorphism
\[
{\rm H}^1(C_H(y),k)  \xleftarrow{{\rm res}^{C_G(x)}_{C_H(x)}}  {\rm H}^1(C_G(y),k)
\]
is zero. This dualises to the required statement.
\end{proof}

Recall that $\mathrm{Core}_G(P)$ and $O_p(G)$, denote the core of $P$ and the $p$-core of G, respectively. If $P$ is a  Sylow $p$-subgroup of $G$, then $O_p(G)=\mathrm{Core}_G(P)$.  We are now able to deduce some further criteria for the nonvanishing of the first Hochschild cohomology of blocks.

\begin{corollary}\label{cor_nonzero_map}
Let $G$ be a finite group and let $B$ be a block of $kG$ with defect group $P$. If for some $x\in P$ the homomorphism ${\rm H}_1(C_P(x),\mathbb{F}_p)\to {\rm H}_1(C_G(x),\mathbb{F}_p)$ is nonzero and for any element $y=hxh^{-1}$ such that $y\notin P$  the homomorphism ${\rm H}_1(C_H(y),\mathbb{F}_p)\to {\rm H}_1(C_G(y),\mathbb{F}_p)$
is zero, then $\HH^1(B,B)\neq0$.
\end{corollary}

\begin{proof}
This follows by combining \cref{prop_centraliser_decomp} with \cref{thm_A_general}.
\end{proof}

\begin{corollary}\label{cor_nonzero_map2}
Let $G$ be a finite group and let $B$ be a block of $kG$ with defect group $P$. If for some element $x \in \mathrm{Core}_G(P)$, the homomorphism ${\rm H}_1(C_P(x),\mathbb{F}_p)\to {\rm H}_1(C_G(x),\mathbb{F}_p)$ is nonzero, then $\HH^1(B)\neq0$.
\end{corollary}

\begin{proof}
This follows from \cref{cor_nonzero_map}.
\end{proof}

\begin{corollary}\label{cor_nonzero_Sylow_old}
Let $G$ be a finite group. If there is a $p$-element $x \in O_p(G)$ such that ${\rm H}_1(C_G(x),\mathbb{F}_p)\neq 0$, then any block $B$ of $kG$ of maximal defect group satisfies $\HH^1(B)\neq 0$.
\end{corollary}

\begin{proof}
Since $x$ is a $p$-element, we may choose a Sylow subgroup $P\subseteq G$ containing $x$, and we may  arrange that $C_P(x)=C_G(x)\cap P$ is a Sylow subgroup of $C_G(x)$. Using the transfer map, this implies that ${\rm H}_1(C_P(x),\mathbb{F}_p)\to {\rm H}_1(C_G(x),\mathbb{F}_p)$ is split surjective, so the result follows from \cref{cor_nonzero_map2}.
\end{proof}

\cref{cor_nonzero_Sylow} in particular holds if there is $x \in O_p(G)$ which  is a strong non-Schur element.

\begin{corollary}
Let $G$ be a finite group and $P$ a non-trivial Sylow $p$-subgroup of $G$. If $O_p(G)\cap Z(P)$ is not contained in $[P,P]$, then any block $B$ of $kG$ of maximal defect group satisfies $\HH^1(B)\neq 0$. 
\end{corollary}
\begin{proof}
    Choose $x \in O_p(G)\cap Z(P)$ which is not in $[P,P]$. Then $P$ is a Sylow $p$-subgroup of $C_G(x)$. 
    Note that $x$ is a strong non-Schur element. Indeed, by the proof of \cref{lem:CPx}, we have $\operatorname{tr}^{C_G(x)}_P(x)= x^{[C_G(x):P]}$ is a non-zero element in $P/[P,P]$ since $p$ does not divide $[C_G(x):P]$. Clearly $x$ is not in $[C_G(x),C_G(x)]$, otherwise the transfer map would be zero.  
\end{proof}

\section{Twisted group algebras}\label{s:tga}

Twisted group algebras are defined by modifying the multiplication in a group algebra using a $2$-cocycle. More precisely, if $\alpha \in Z^2(G; k^\times)$ is a $2$-cocycle on $G$, then $k_\alpha G$ is the algebra with a $k$-basis $\{\hat{g} \mid g \in G\}$
and with multiplication determined by
\[
\hat{g} \cdot \hat{h} = \alpha(g, h) \widehat{gh} \quad \text{for } g, h \in G.
\]

There are some families of blocks which are Morita equivalent to twisted group algebras.  For instance, if $B$ is a  block with  normal defect group in $G$, then $B$ is Morita equivalent to a twisted group algebra \cite{Ku}. This is illustrated by the following example:

A family of subgroups of $G$ has the trivial intersection (TI) property if any two members $P, Q$ satisfy $P\cap Q=\{1\}$. If a block $B$ of $kG$ has TI defect groups, then there exists a stable equivalence of Morita type between $B$ and its Brauer correspondent $b$ in $kN_G(P)$, for any defect group $P$ \cite[\S   11.2]{KZ98}. Since $b$ has defect group $P$ normal in $N_G(P)$ we have that $b$ is Morita equivalent to a twisted group algebra. The invariance of Hochschild cohomology under stable equivalences of Morita type together with \cite[Proposition 4.1]{Todea} implies $\HH^1(B)\neq 0$. 

For these families of blocks, Linckelmann's conjecture reduces to the following question:

\begin{question} \cite[Question 1.3]{Todea}\label{Q:Todea}
    Is the first Hochschild cohomology of twisted group algebras of finite groups
(with order divisible by $p$) nontrivial? If not, can someone find examples of such twisted
group algebras which have zero first Hochschild cohomology?
\end{question}

A partial affirmative answer to this question has been established in several cases, including 
p-solvable groups \cite{Todea}, all finite simple groups \cite{Murphy}, and certain dihedral groups \cite{OteroSanchez25}. 
\begin{definition}
    An element $x$ of $G$ is called \emph{$\alpha$-regular} if $\alpha(g,x)=\alpha(x,g)$ for all $g\in C_G(x)$. 
\end{definition}

For such an element one can identify ${\rm H}^1(C_G(x),k_\alpha \hat{x})={\rm H}^1(C_G(x),k)$.

\begin{lemma}\cite[Lemma 2.11]{Murphy}\label{lem:tga}
If $G$ is a finite group and $\alpha \in Z^2(G; k^\times)$, then every non-Schur $p$-element $x$ is $\alpha$-regular. In particular, if $G$ contains a  non-Schur $p$-element  then $\mathrm{HH}^1(k_\alpha G,k_\alpha G) \neq 0$.
\end{lemma}

As a consequence of Theorem \ref{thm:g7} and Lemma \ref{lem:tga} observe that Question \ref{Q:Todea} has a positive answer in all characteristics greater than $5$.

\begin{corollary}\label{thm:gt7}
    If $G$ is a finite group and $k$ is a field of prime characteristic greater than $5$, then  $\mathrm{HH}^1(k_\alpha G,k_\alpha G) \neq 0$ for any $2$-cocycle $\alpha$ on $G$. \qed
\end{corollary}

The following result is a generalisation of \cite[Proposition 4.1 (d)]{Todea}; the terminology used is explained in \cite[\S 8.2]{Linckelmann18-II}.

\begin{corollary}\label{cor_F_centric}
Let $G$ be a finite group and $k$ is a field of prime characteristic $p>5$. Let $B=kGb$ be a block of $kG$ with nontrivial defect group $P$. Let $Q$ be an $\mathcal{F}$-centric subgroup of $P$ and let $f$ be the block-idempotent of $kC_G(Q)$ such that $(Q, f ) \leq (P, b)$. Assume that $G = N_G(Q,f)$. Then $\HH^1(B)\neq 0$.
\end{corollary}
  
\begin{proof}
  It follows from \cite[Corollary~8.12.9]{Linckelmann18-II} that  $B$ is Morita equivalent to a twisted group algebra $k_{\alpha}L$, where $L$ is a finite group having $P$ as a Sylow $p$-subgroup. The statement follows from
  Corollary \ref{thm:gt7}.
\end{proof}

Following \cite{FJL}, an element $x
\in G$ has the \emph{commutator index property} if $p$ divides the index of $[C_G(x),C_G(x)]$ in $C_G(x)$, or equivalently, if ${\rm H}_1(C_G(x),\mathbb{F}_p)\neq 0$.

\begin{remark}
 As Todea pointed out in \cite[Page 3]{Todea}, if $G$ possess an element with the commutator index property 
 then $\mathrm{HH}^1(k_\alpha G,k_\alpha G) \neq 0$ but  the converse does not hold. For example, if $G$ is the Rudvalis group $Ru$, then by \cite[Proposition 3.3]{Murphy} we have that  $\mathrm{HH}^1(k_\alpha G,k_\alpha G) \neq 0$. However, using GAP one can check  that   $G=Ru$ does not admit any element with the commutator index property for $p=3$, see also \cite[Page 3]{Todea}. 
\end{remark}

\bibliographystyle{alpha}
\bibliography{references}

\end{document}